\documentclass{amsart}
\usepackage{amssymb}
\usepackage{mathtools}
\usepackage{tikz-cd}
\DeclareMathOperator{\Vect}{Vect}
\newcommand{\id}{\mathrm{id}}
\DeclareMathOperator{\lproj}{lproj}
\DeclareMathOperator{\rproj}{rproj}
\newcommand{\btensor}{\mathbin{\bar{\otimes}}}

\newtheorem{theorem}{Theorem}[section]
\newtheorem{lemma}[theorem]{Lemma}
\newtheorem{proposition}[theorem]{Proposition}
\newtheorem{corollary}[theorem]{Corollary}

\theoremstyle{definition}
\newtheorem{definition}[theorem]{Definition}

\theoremstyle{remark}
\newtheorem{examples}[theorem]{Examples}
\newtheorem{example}[theorem]{Example}
\newtheorem{remark}[theorem]{Remark}

\numberwithin{equation}{section}

\title[Lax Centers and Monoidal Functors]{Lax Centers and Monoidal Functors for Comodule Categories Associated with Acyclic Quivers}



\author{Wentao Zhao}
\author{Jingheng Zhou}

\address{Wentao Zhao: School of Mathematics, 
Shanghai University of Finance and Economics, Shanghai 200433, 
Shanghai, China}

\email{2025213182@stu.sufe.edu.cn}

\address{Jingheng Zhou: School of Mathematics, 
Shanghai University of Finance and Economics, Shanghai 200433, 
Shanghai, China}

\email{zhoujingheng@mail.shufe.edu.cn}

\begin{document}

\begin{abstract}
In this paper, we study the monoidal category of finite-dimensional right comodules over the path algebra $\Bbbk Q$ of a finite acyclic quiver $Q$, equipped with the weak bialgebra structure in which every path is grouplike. We prove that the lax center of this category is monoidally equivalent to $\operatorname{Rep}^{\mathrm{fd}}_{\Bbbk}(Q)$, the category of finite-dimensional representations of $Q$. We also classify all $\Bbbk$-linear strong monoidal functors between any two such comodule categories, and all monoidal natural transformations between them. For the right adjoint of such a functor, we determine when
the left or right projection formula holds.
When both projection formulas hold, this adjoint will induce a functor between the lax centers,
 and we describe the induced functor explicitly on quiver representations.
\end{abstract}

\maketitle

\section{Introduction}

Weak bialgebras and weak Hopf algebras were introduced by B\"ohm, Nill, and
Szlach\'anyi in 1999 to provide a common framework for groupoids, face algebras,
and quantum groupoids. In a weak bialgebra the comultiplication need not be
unital and the counit need not be multiplicative, but the category of right
comodules is still monoidal \cite{bohm2011weak}. This makes weak bialgebras a natural
source of monoidal categories that are not rigid, and it motivates the study of
monoidal centers in this setting.

The theory of monoidal centers has a long history. The Drinfeld center of a
monoidal category was developed in the theory of braided tensor categories and
has been studied from several points of view
\cite{JS1993,S2001,S2004}. For a rigid monoidal category, the Drinfeld center
is a braided monoidal category and is a fundamental invariant of the monoidal
structure. In 2005, Caenepeel, Wang, and Yin studied the weak center of modules
over a weak bialgebra and related it to Yetter--Drinfeld modules
\cite{CWY2005}. In 2007, Day, Panchadcharam, and Street studied the lax center,
in which the half-braiding is not required to be invertible \cite{DPS2007}. For
rigid categories the lax center coincides with the Drinfeld center, but for
non-rigid categories the two can differ substantially.

Quiver representations provide a concrete and important class of examples. For
a finite acyclic quiver $Q$, the path algebra $\Bbbk Q$ has a basis indexed by
paths and carries the canonical weak bialgebra structure in which every path is
grouplike \cite[Lemma~2.1(1)]{ZZ2022}. Its finite-dimensional comodules are
precisely the path-graded vector spaces, so the comodule category
$\mathcal C_Q$ is finite semisimple with simple objects indexed by paths. The
tensor product is governed by concatenation of paths. Tensor structures on
quiver representations were studied by Herschend in 2008 and by Huang and
Torrecillas in 2013 \cite{H2008,HT2013}, and the canonical weak bialgebra
structure was used by Zhang and Zhou in 2022 in their Frobenius--Perron theory
of quiver representations \cite{ZZ2022}. However, when $Q$ has an arrow,
$\mathcal C_Q$ is not rigid: its unit object is a direct sum of simple objects
indexed by trivial paths, while tensoring the simple object indexed by an arrow
with any object cannot produce such a summand. The usual identification of the
lax center with the Drinfeld center is therefore unavailable, and the lax center
of $\mathcal C_Q$ requires a separate study.

More recently, Flake, Laugwitz, and Posur developed a general theory of
projection formula morphisms for right adjoints of strong monoidal functors and
described when these adjoints induce functors between Drinfeld centers
\cite{flake2026projection}. They also studied Frobenius monoidal structures arising from
ambiadjunctions and their lifts to centers \cite{FLP2025}. These results provide
the tools for understanding the functorial behavior of monoidal centers, but
they do not by themselves describe the lax center of $\mathcal C_Q$ or classify
monoidal functors between the categories $\mathcal C_Q$ and $\mathcal C_R$. The
present paper addresses these questions.

In this paper, we study the lax center of $\mathcal C_Q$ and the structure of
strong monoidal functors between such comodule categories. We find that there is a $\Bbbk$-linear monoidal equivalence between the lax center of $\mathcal C_Q$ and 
the representation category of $Q$:
\[
Z_{\mathrm{lax}}(\mathcal C_Q)
\simeq
\operatorname{Rep}^{\mathrm{fd}}_{\Bbbk}(Q),
\]
Also, the equivalence carries the lax braiding to the pointwise symmetric braiding and 
equivalence can restrict to a braided monoidal equivalence
$Z(\mathcal C_Q)\simeq
\operatorname{Rep}^{\mathrm{fd}}_{\Bbbk,\mathrm{iso}}(Q)$.
As for two quivers $Q$ and $R$, we classifies all $\Bbbk$-linear strong monoidal functors
$\mathcal C_Q\to\mathcal C_R$ in terms of decompositions of the vertex set
$R_0$ indexed by $Q_0$ and objects $X_a\in\mathcal C_R$ associated with the
arrows of $Q$. We also describe all monoidal natural transformations between
such functors. Last, we give necessary and sufficient conditions
for the left and right projection formulas of the right adjoints of these
functors to be isomorphisms. When both formulas hold, we obtain a functor
between the corresponding lax centers and describe the induced functor
explicitly on quiver representations.

The paper is organized as follows. Section~2 recalls the weak bialgebra
structure on $\Bbbk Q$ and describes $\mathcal C_Q$ as a category of
path-graded vector spaces. Section~3 computes the lax center and the Drinfeld
center of $\mathcal C_Q$. Section~4 classifies strong monoidal functors between
$\mathcal C_Q$ and $\mathcal C_R$, together with their monoidal natural
transformations. Section~5 studies the projection formulas for the right
adjoints of these functors and describes the induced functors on centers.

\section{Preliminaries}

In this section, we recall the weak bialgebra structure on
$\Bbbk Q$ and describe its comodule category.

Our conventions for monoidal categories, strong monoidal functors, and
coherence follow
\cite{ML1998,EGNO2015}.

We use the standard definition of a weak bialgebra over $\Bbbk$, see
\cite{BNS1999}.

Let
$Q=(Q_0,Q_1,s,t)$
be a finite quiver, where $Q_0$ is the set of vertices and $Q_1$ is the set of arrows.
For each $a\in Q_1$, $s(a)$ and $t(a)$ denote the source and target of $a$, respectively.
All quivers considered in this paper are assumed to have nonempty vertex sets.
Let $\mathcal P(Q)$ denote the set of all paths in $Q$.
For each vertex $i\in Q_0$, let $e_i$ denote the trivial path of length $0$ at $i$.
For any path $p$, let $s(p)$ and $t(p)$ denote its source and target, respectively,
and let $\ell(p)$ denote its length.
A path of positive length is called nontrivial.

As a vector space, the path algebra of $Q$ is
$\Bbbk Q=\bigoplus_{p\in\mathcal P(Q)}\Bbbk p$.
For $p,q\in\mathcal P(Q)$, define their product by
\[
pq=
\begin{cases}
\text{the path obtained by first traversing $p$ and then $q$},
&t(p)=s(q),\\
0,&t(p)\neq s(q).
\end{cases}
\]
Extending this product bilinearly makes $\Bbbk Q$ an associative
$\Bbbk$-algebra
\cite{ASS2006,DW2005,
S2014}.
When $Q$ is acyclic and $Q_0$ is finite, every path has length at most
$|Q_0|-1$, so $\mathcal P(Q)$ is finite and $\Bbbk Q$ is a
finite-dimensional algebra with unit
\[
1_{\Bbbk Q}=\sum_{i\in Q_0}e_i.
\]
Denote this multiplication by $m$ and the unit map by
$\eta\colon \Bbbk\to \Bbbk Q$.

Define $\Delta(p)=p\otimes p$ and $\varepsilon(p)=1$ for
$p\in\mathcal P(Q)$, and extend these maps by $\Bbbk$-linearity.
It is known that the path algebra $\Bbbk Q$ with the above coalgebra structure is a cocommutative weak bialgebra
\cite[Lemma~2.1(1)]{ZZ2022}.

Fix a finite acyclic quiver $Q$, and write
$\mathcal C_Q\coloneqq\operatorname{Comod}^{\mathrm{fd}}(H(Q))$.
We write $\btensor$ for the tensor product in $\mathcal C_Q$,
reserving $\otimes_{\Bbbk}$ for the usual tensor product of $\Bbbk$-vector spaces.

Let $S_p=\Bbbk u_p$ be the one-dimensional right comodule satisfying
$\rho_{S_p}(u_p)=u_p\otimes_{\Bbbk} p$. For $(M,\rho_M)\in\mathcal C_Q$ and $p\in\mathcal P(Q)$, set
\[
M_p\coloneqq\{m\in M\mid\rho_M(m)=m\otimes_{\Bbbk} p\}.
\]
We call $M_p$ the homogeneous component of $M$ indexed by $p$.

For a finite-dimensional $\Bbbk$-vector space $V$ and $X\in\mathcal C_Q$, we write
$V\otimes_{\Bbbk}X$ for the comodule with coaction
$\id_V\otimes_{\Bbbk}\rho_X$, so that $(V\otimes_{\Bbbk}X)_p=V\otimes_{\Bbbk}X_p$ for all
$p\in\mathcal P(Q)$.
We use the same convention for the comodule categories associated with other
quivers.

Let $\Vect_{\mathcal P(Q)}^{\mathrm{fd}}$ denote the category of finite-dimensional
vector spaces graded by the set $\mathcal P(Q)$ of paths. Next proposition tells us 
$\mathcal C_Q$ is equivalent to $\Vect_{\mathcal P(Q)}^{\mathrm{fd}}$.

\begin{proposition}
\label{xxprop2.1}
Every $M\in\mathcal C_Q$ has a natural finite direct sum decomposition
$M=\bigoplus_{p\in\mathcal P(Q)}M_p$,
and the comodule morphisms are precisely the linear maps that preserve the homogeneous components indexed by the paths.
The standard monoidal structure on the category of right $H(Q)$-comodules satisfies
\begin{equation}\label{xxeq2-1}
(M\btensor N)_z
=
\bigoplus_{\substack{p,q\in\mathcal P(Q)\\t(p)=s(q),\ pq=z}}
M_p\otimes_{\Bbbk}N_q,
\qquad \text{ and }
\mathbf1_Q=\bigoplus_{v\in Q_0}S_{e_v},
\end{equation}
where $\mathbf1_Q$ is the unit object of $\mathcal C_Q$.

In particular,
\[
S_p\btensor S_q
\cong
\begin{cases}
S_{pq},&t(p)=s(q),\\
0,&t(p)\neq s(q).
\end{cases}
\]

Moreover, $\mathcal C_Q$ is a finite semisimple $\Bbbk$-linear category,
$\{S_p\mid p\in\mathcal P(Q)\}$ is a set of representatives for the isomorphism classes of simple objects, and
$\mathcal C_Q\simeq\Vect_{\mathcal P(Q)}^{\mathrm{fd}}$
as monoidal categories, where the tensor product on the right is given by
\eqref{xxeq2-1}.
\end{proposition}

\begin{proof}
Write $\rho_M(m)=\sum_p m_p\otimes_{\Bbbk} p$
with respect to the path basis.
The counit axiom and coassociativity together with
$\Delta(p)=p\otimes p$ give, respectively,
\[
m=\sum_p m_p,
\qquad
\rho_M(m_p)=m_p\otimes_{\Bbbk} p.
\]

The linear independence of the paths shows that the sum is direct.
The same computation shows that comodule morphisms are precisely the
linear maps preserving each $M_p$. Conversely, for a graded vector space
$V=\bigoplus_pV_p$, define
$\rho(v)=v\otimes_{\Bbbk} p$ for $v\in V_p$.
This gives a quasi-inverse and hence a $\Bbbk$-linear equivalence of categories.

In the monoidal structure of
\cite[Section~3]{bohm2011weak}, for $m\in M_p$ and
$n\in N_q$, the idempotent defining $M\btensor N$ becomes
$P_{M,N}(m\otimes_{\Bbbk} n)=\varepsilon(pq)m\otimes_{\Bbbk} n$.
It is the identity when $t(p)=s(q)$ and is zero otherwise, giving the stated formula for the tensor product.
The unit object arising from the same general construction becomes
$\bigoplus_{v\in Q_0}\Bbbk e_v\cong\bigoplus_vS_{e_v}$.

The formula for tensor products of simple objects follows immediately.
Every object is also naturally isomorphic to
$\bigoplus_pM_p\otimes_{\Bbbk}S_p$,
which shows that $\mathcal C_Q$ is finite semisimple and gives the
description of its simple objects.

Finally, the above equivalence sends the tensor product and the unit object,
respectively, to the tensor product described by
\eqref{xxeq2-1}
and to the unit object whose component is $\Bbbk$ at each trivial path and zero
elsewhere.
The associativity and unit constraints are the canonical vector-space isomorphisms on each homogeneous component.
Thus the equivalence is a monoidal equivalence.
\end{proof}

For $B\subseteq Q_0$, set
$\mathbf1_B\coloneqq\bigoplus_{v\in B}S_{e_v}$, and let
\[
\iota_B\colon\mathbf1_B\longrightarrow\mathbf1_Q,
\qquad
\pi_B\colon\mathbf1_Q\longrightarrow\mathbf1_B,
\qquad
E_B\coloneqq\iota_B\pi_B
\]
be the canonical inclusion, the canonical projection, and the corresponding
idempotent in $\operatorname{End}_{\mathcal C_Q}(\mathbf1_Q)$.
We write $\iota_v$, $\pi_v$ and $E_v$ when $B=\{v\}$, and add a superscript
indicating the quiver when needed, as in $\mathbf1_B^R$, $\iota_B^R$,
$\pi_B^R$ and $E_B^R$.
Under
$\operatorname{End}_{\mathcal C_Q}(\mathbf1_Q)\cong \Bbbk^{Q_0}$, the idempotent
$E_B$ corresponds to the indicator vector of $B$. In particular,
$B\mapsto E_B$ is a bijection from the subsets of $Q_0$ onto the idempotents
of $\operatorname{End}_{\mathcal C_Q}(\mathbf1_Q)$.

\section{The Lax Center and Quiver Representations}

\begin{definition}\label{xxdef3.1}
Let $\mathcal C$ be a monoidal category.
We denote the left and right unit constraints by
$l_X\colon\mathbf1\otimes X\to X$ and
$r_X\colon X\otimes\mathbf1\to X$, respectively.
An object of the lax center
$Z_{\mathrm{lax}}(\mathcal C)$
is a pair $(V,\gamma)$, where $V\in\mathcal C$ and
$\gamma\colon V\otimes-\Rightarrow-\otimes V$
is a natural transformation, called a lax half-braiding, whose component at
$X\in\mathcal C$ is
$\gamma_X\colon V\otimes X\longrightarrow X\otimes V$.
The natural transformation $\gamma$ is required to satisfy the unit condition
and, suppressing the associativity constraints in the notation, the
compatibility condition
\[
\begin{aligned}
\gamma_{\mathbf1}
&=l_V^{-1}\circ r_V,\\
\gamma_{X\otimes Y}
&=(\id_X\otimes\gamma_Y)
  \circ(\gamma_X\otimes\id_Y).
\end{aligned}
\]
The morphisms $\gamma_X$ are not required to be invertible.
A morphism from $(V,\gamma)$ to $(W,\delta)$ is a morphism
$f\colon V\to W$ such that
$(\id_X\otimes f)\circ\gamma_X
=
\delta_X\circ(f\otimes\id_X)$
for all~$X$.
The objects of the lax center for which all $\gamma_X$ are invertible form
the Drinfeld center $Z(\mathcal C)$.
\end{definition}

We use the monoidal structure on the lax center described in
\cite{DPS2007}.
Suppressing associativity constraints, its tensor product is
\[
(V,\gamma)\btensor(W,\delta)
=
\bigl(V\btensor W,
(\gamma_X\btensor\id_W)\circ(\id_V\btensor\delta_X)\bigr),
\]
and its unit object is
$\bigl(\mathbf1_Q,(r_X^{-1}\circ l_X)_X\bigr)$.

Let $\operatorname{Rep}^{\mathrm{fd}}_{\Bbbk}(Q)$ denote the category of finite-dimensional
representations of $Q$ over $\Bbbk$.
Let $\operatorname{Rep}^{\mathrm{fd}}_{\Bbbk,\mathrm{iso}}(Q)$ denote the full
subcategory in which every arrow map is an isomorphism.
We write objects of $\operatorname{Rep}^{\mathrm{fd}}_{\Bbbk}(Q)$ as
$\bigl((V_v)_{v\in Q_0},
(\varphi_a\colon V_{s(a)}\to V_{t(a)})_{a\in Q_1}\bigr)$.

This category is equipped with the pointwise tensor product
\[
(V\otimes W)_v
=
V_v\otimes_{\Bbbk}W_v,
\qquad
(\varphi^{V\otimes W})_a
=
\varphi^V_a\otimes_{\Bbbk}\varphi^W_a.
\]
Its unit object takes the value $\Bbbk$ at each vertex and the identity map at
each arrow.

The lax center is lax braided, with
$c_{(V,\gamma),(W,\delta)}\coloneqq\gamma_W$
\cite{DPS2007}.
On $\operatorname{Rep}^{\mathrm{fd}}_{\Bbbk}(Q)$, we use the pointwise
symmetric braiding whose component at a vertex $v$ is the flip
$V_v\otimes_{\Bbbk} W_v\to W_v\otimes_{\Bbbk} V_v$.

\begin{theorem}\label{xxthm3.2}
Let $Q$ be a finite acyclic quiver.
There is a $\Bbbk$-linear monoidal equivalence
\[
\Phi_Q\colon
Z_{\mathrm{lax}}(\mathcal C_Q)
\xrightarrow{\sim}
\operatorname{Rep}^{\mathrm{fd}}_{\Bbbk}(Q)
\]
that carries the lax braiding of $Z_{\mathrm{lax}}(\mathcal C_Q)$ to the
pointwise symmetric braiding. In particular, this lax braiding is
invertible. Moreover, $\Phi_Q$ restricts to a braided monoidal equivalence
$Z(\mathcal C_Q)\xrightarrow{\sim}
\operatorname{Rep}^{\mathrm{fd}}_{\Bbbk,\mathrm{iso}}(Q)$.
\end{theorem}

\begin{proof}
Let
$(V,\gamma)\in Z_{\mathrm{lax}}(\mathcal C_Q)$,
and write
$V=\bigoplus_{p\in\mathcal P(Q)}V_p$.

For each $v\in Q_0$, naturality of $\gamma$ with respect to
$E_v\colon\mathbf1_Q\to\mathbf1_Q$ is the commutativity of the square
\[
\begin{tikzcd}[column sep=large,row sep=large]
V\btensor\mathbf1_Q
  \arrow[r,"\gamma_{\mathbf1_Q}"]
  \arrow[d,"\id_V\btensor E_v"']
&
\mathbf1_Q\btensor V
  \arrow[d,"E_v\btensor\id_V"]\\
V\btensor\mathbf1_Q
  \arrow[r,"\gamma_{\mathbf1_Q}"']
&
\mathbf1_Q\btensor V.
\end{tikzcd}
\]
Since the unit condition gives
$\gamma_{\mathbf1_Q}=l_V^{-1}\circ r_V$, substituting this equality into
the square, composing with $l_V$ on the left and with $r_V^{-1}$ on the
right, we obtain
\[
l_V(E_v\btensor\id_V)l_V^{-1}
=
r_V(\id_V\btensor E_v)r_V^{-1}.
\]
Suppose
$p\colon i\to j$
and
$V_p\neq0$.
On $V_p$, the left-hand side is the identity if and only if $v=i$, whereas the right-hand side is the identity if and only if $v=j$.
Thus $i=j$.
Since $Q$ has no oriented cycles, such a path must be trivial.
Therefore $V=\bigoplus_{v\in Q_0}V_{e_v}$.

For a path
$p\colon i\to j$,
Proposition~\ref{xxprop2.1} gives
\[
V\btensor S_p
\cong
V_{e_i}\otimes_{\Bbbk}S_p,
\qquad
S_p\btensor V
\cong
S_p\otimes_{\Bbbk}V_{e_j}.
\]
Thus $\gamma_{S_p}$ corresponds uniquely to a linear map
$\varphi_p\colon V_{e_i}\to V_{e_j}$.

The unit and compatibility conditions in
Definition~\ref{xxdef3.1} give
$\varphi_{e_v}=\id_{V_{e_v}}$ and
$\varphi_{pq}=\varphi_q\circ\varphi_p$ for all composable paths $p,q$.
Hence each $\varphi_p$ is the composite of the maps $\varphi_a$ along $p$,
and we associate to $(V,\gamma)$ the finite-dimensional representation
$\bigl((V_{e_v})_v,(\varphi_a)_a\bigr)$ of $Q$.

Conversely, given a representation
$\bigl((V_v)_v,(\varphi_a)_a\bigr)$, let
$V=\bigoplus_vV_v$,
with $V_v$ placed in the homogeneous component indexed by the trivial
path $e_v$.
For each path $p$, let $\varphi_p$ denote the composite of the arrow maps along
$p$.
For $M\in\mathcal C_Q$ and $x\in M_p$, define
\[
\gamma_M(y\otimes_{\Bbbk} x)
=
x\otimes_{\Bbbk}\varphi_p(y)
\qquad
(y\in V_{e_{s(p)}}).
\]

These maps preserve the $\mathcal P(Q)$-grading and are natural in $M$.
The identities $\varphi_{e_v}=\id_{V_v}$ and
$\varphi_{pq}=\varphi_q\circ\varphi_p$ give the unit and compatibility
conditions, respectively, in Definition~\ref{xxdef3.1}.
They therefore define an object of the lax center.

Now let $(V,\gamma)$ be an object of the lax center, and let
$\varphi_p$ be the linear maps obtained from it as above.
For $M\in\mathcal C_Q$ and $x\in M_p$, let
$f_x\colon S_p\to M$ be defined by $f_x(u_p)=x$.

Naturality of $\gamma$ with respect to $f_x$ shows that $\gamma_M$ is
given by the defining formula above, so the reverse construction
recovers the original $\gamma$.
The two constructions are also mutually inverse on morphisms: morphisms
in the lax center correspond precisely to morphisms of representations.
This gives a $\Bbbk$-linear equivalence of categories.

Under this equivalence, the tensor product of
$(V,\gamma)$ and $(W,\delta)$ has vertex components
$V_v\otimes_{\Bbbk}W_v$, and the map associated with an arrow $a$ is
$\varphi^V_a\otimes_{\Bbbk}\varphi^W_a$.
The unit object corresponds to the unit object of
$\operatorname{Rep}^{\mathrm{fd}}_{\Bbbk}(Q)$. Hence the equivalence is monoidal.
Since $W$ is supported only on trivial paths, the lax braiding
$c_{(V,\gamma),(W,\delta)}=\gamma_W$ satisfies
$\gamma_W(y\otimes_{\Bbbk} x)=x\otimes_{\Bbbk}\varphi_{e_v}(y)=x\otimes_{\Bbbk} y$
for $y\in V_{e_v}$ and $x\in W_{e_v}$.
Thus the lax braiding is invertible and corresponds to the
pointwise symmetric braiding on
$\operatorname{Rep}^{\mathrm{fd}}_{\Bbbk}(Q)$.

Finally,
$\gamma_{S_a}$
is invertible if and only if
$\varphi_a$
is invertible.
Together with the formula for composition along paths, this shows that
all arrow maps are invertible if and only if all
$\gamma_X$
are invertible.
Hence the above monoidal equivalence restricts to
$Z(\mathcal C_Q)\simeq
\operatorname{Rep}^{\mathrm{fd}}_{\Bbbk,\mathrm{iso}}(Q)$.
\end{proof}

\begin{remark}\label{xxrem3.3}

Acyclicity enters
Theorem~\ref{xxthm3.2} in two places. First,
it makes $\mathcal P(Q)$ finite, so that $\mathcal C_Q$ is finite
semisimple (Proposition~\ref{xxprop2.1}). Second, in
the proof, for a nonzero component $V_p$ with $p\colon i\to j$, the
identity
\[
l_V(E_v\btensor\id_V)l_V^{-1}
=
r_V(\id_V\btensor E_v)r_V^{-1}
\]
forces $i=j$. Since a nontrivial path from a vertex to itself is an
oriented cycle, acyclicity implies that $p$ is trivial.

If $Q$ is allowed to have oriented cycles, then $\mathcal P(Q)$ is
infinite, $\mathcal C_Q$ is not finite semisimple, and objects of the lax
center may be supported on nontrivial paths. Consequently the conclusion
$V=\bigoplus_{v\in Q_0}V_{e_v}$, and hence the description of
Theorem~\ref{xxthm3.2}, fails in general.
Thus acyclicity is essential for this description, not merely a technical
convenience.

\end{remark}

\section{Classification of Strong Monoidal Functors}

For the remainder of the paper, let $Q$ and $R$ be finite acyclic quivers,
and write $\mathcal C_R\coloneqq\operatorname{Comod}^{\mathrm{fd}}(H(R))$.
By a decomposition $R_0=\bigsqcup_{i\in Q_0}B_i$ we mean a family of
pairwise disjoint, possibly empty subsets $B_i\subseteq R_0$ whose union is
$R_0$.
For a strong monoidal functor $F$, we write
$J_{M,N}\colon F(M)\btensor F(N)\xrightarrow{\sim}F(M\btensor N)$ and
$J_0\colon\mathbf1_R\xrightarrow{\sim}F(\mathbf1_Q)$ for its tensor and unit
isomorphisms, respectively.

For $r\in\mathcal P(R)$, write $S_r^R$ for the corresponding
one-dimensional simple object of $\mathcal C_R$.

For any family of objects $\mathcal A$ in $\mathcal C_R$, let
$\operatorname{add}(\mathcal A)$
denote the full subcategory whose objects are isomorphic to direct summands of finite direct sums of objects in $\mathcal A$.
For
$B,C\subseteq R_0$,
set
\[
\mathcal C_R[B,C]
\coloneqq
\operatorname{add}
\left\{
S_r^R\ \middle|\
r\in\mathcal P(R),\
s(r)\in B,\
t(r)\in C
\right\}.
\]
If $B$ or $C$ is empty, then $\mathcal C_R[B,C]$
consists only of zero objects.

Let
$\operatorname{Fun}^{\otimes}_{\Bbbk}(\mathcal C_Q,\mathcal C_R)$
denote the category whose objects are $\Bbbk$-linear strong monoidal functors
and whose morphisms are monoidal natural transformations.
We write $\operatorname{Hom}_{\otimes}(F,G)$ for the set of morphisms from
$F$ to $G$ in this category.
For such an $F$, set
$\operatorname{End}_{\otimes}(F)\coloneqq\operatorname{Hom}_{\otimes}(F,F)$,
and let $\operatorname{Aut}_{\otimes}(F)$ denote the group of units of this
monoid.

\begin{theorem}\label{xxthm4.1}
Let $Q,R$ be finite acyclic quivers.
There is an equivalence of categories
\[
\operatorname{Fun}^{\otimes}_{\Bbbk}(\mathcal C_Q,\mathcal C_R)
\simeq
\coprod_{R_0=\bigsqcup_{i\in Q_0}B_i}
\prod_{a\in Q_1}
\mathcal C_R[B_{s(a)},B_{t(a)}],
\]
under which $((B_i),(X_a))$ corresponds to the functor $F_{B,X}$
constructed in Lemma~\ref{xxlem4.3} below.
On the right-hand side, a morphism
$((B_i),(X_a))\to((B_i'),(Y_a))$ exists only when $B_i=B_i'$ for all $i$,
and is then a family $(f_a\colon X_a\to Y_a)_{a\in Q_1}$.
We regard the strong monoidal structure as part of the data of an object of
$\operatorname{Fun}^{\otimes}_{\Bbbk}(\mathcal C_Q,\mathcal C_R)$.
\end{theorem}

\begin{lemma}\label{xxlem4.2}
Let
$(F,J,J_0)\colon\mathcal C_Q\longrightarrow\mathcal C_R$
be a $\Bbbk$-linear strong monoidal functor.
Then there is a unique decomposition $R_0=\bigsqcup_{i\in Q_0}B_i$
such that, for every $i\in Q_0$,
\[
F(S_{e_i})
\cong
A_i
\coloneqq
\bigoplus_{v\in B_i}S_{e_v}^R.
\]

Moreover,
$F(S_a)\in\mathcal C_R[B_i,B_j]$
for each arrow $a\colon i\to j$.

\end{lemma}

\begin{proof}
For each $i\in Q_0$, the idempotent $J_0^{-1}F(E_i^Q)J_0$ of
$\operatorname{End}_{\mathcal C_R}(\mathbf1_R)$ equals $E_{B_i}^R$ for a
unique subset $B_i\subseteq R_0$. Since the $E_i^Q$ are pairwise orthogonal
with sum $\id_{\mathbf1_Q}$, the $E_{B_i}^R$ are pairwise orthogonal with sum
$\id_{\mathbf1_R}$, so the $B_i$ form the stated decomposition.

Note that $A_i=\mathbf1_{B_i}^R$. Then
\[
\theta_{e_i}\coloneqq F(\pi_i^Q)J_0\iota_{B_i}^R
\colon A_i\xrightarrow{\sim}F(S_{e_i}),
\qquad
\theta_{e_i}^{-1}=\pi_{B_i}^R J_0^{-1}F(\iota_i^Q),
\]
as follows from $\pi_i^Q\iota_i^Q=\id_{S_{e_i}}$,
$\pi_{B_i}^R\iota_{B_i}^R=\id_{A_i}$,
$\iota_{B_i}^R\pi_{B_i}^R=E_{B_i}^R=J_0^{-1}F(E_i^Q)J_0$.
Moreover, $J_0\iota_{B_i}^R=F(\iota_i^Q)\theta_{e_i}$.

For an arrow $a\colon i\to j$, set $X_a\coloneqq F(S_a)$.
From
$S_{e_i}\btensor S_a\cong S_a\cong S_a\btensor S_{e_j}$
and the strong monoidal structure, we obtain
$A_i\btensor X_a\cong X_a$ and
$X_a\btensor A_j\cong X_a$.

Decomposing $X_a$ into simple objects indexed by paths and comparing
multiplicities using
Proposition~\ref{xxprop2.1},
we see that every simple direct summand $S_r^R$ of $X_a$ satisfies
$s(r)\in B_i$ and $t(r)\in B_j$.
Hence $X_a\in\mathcal C_R[B_i,B_j]$.
\end{proof}

\begin{lemma}\label{xxlem4.3}
Given a decomposition $R_0=\bigsqcup_iB_i$ and objects
$X_a\in\mathcal C_R[B_i,B_j]$ $(a\colon i\to j)$, set
$X_{e_i}\coloneqq\bigoplus_{v\in B_i}S_{e_v}^R$.
For a path $p=a_1\cdots a_m$ of positive length, set
$X_p\coloneqq X_{a_1}\btensor\cdots\btensor X_{a_m}
\in\mathcal C_R[B_{s(p)},B_{t(p)}]$,
using a fixed left-associated bracketing.
For an object $M\in\mathcal C_Q$ and a morphism
$f\colon M\to N$ in $\mathcal C_Q$, define
\[
F_{B,X}(M)
\coloneqq
\bigoplus_{p\in\mathcal P(Q)}M_p\otimes_{\Bbbk}X_p,
\qquad
F_{B,X}(f)
\coloneqq
\bigoplus_p
\bigl(f|_{M_p}\otimes_{\Bbbk}\id_{X_p}\bigr).
\]
Then $F_{B,X}\colon\mathcal C_Q\to\mathcal C_R$ admits the structure of a
$\Bbbk$-linear strong monoidal functor satisfying
$F_{B,X}(S_{e_i})\cong X_{e_i}$ and
$F_{B,X}(S_a)\cong X_a$.

Moreover, every $\Bbbk$-linear strong monoidal functor
$F\colon\mathcal C_Q\to\mathcal C_R$
is monoidally isomorphic to the functor
$F_{B,X}$
constructed from the data obtained in
Lemma~\ref{xxlem4.2}.
\end{lemma}

\begin{proof}
Set $A_i\coloneqq X_{e_i}$.
For composable paths $p,q$ of positive length, Mac Lane's coherence theorem
\cite{ML1998}
gives a canonical isomorphism
$\kappa_{p,q}\colon X_p\btensor X_q\xrightarrow{\sim}X_{pq}$.

If
$p\colon i\to j$,
define the isomorphisms involving trivial paths by
\[
\kappa_{e_i,p}
\coloneqq
l_{X_p}
(\iota_{B_i}^R\btensor\id_{X_p}),
\qquad
\kappa_{p,e_j}
\coloneqq
r_{X_p}
(\id_{X_p}\btensor\iota_{B_j}^R).
\]
Since $A_h\btensor X_p=0$ for $h\neq i$ and $X_p\btensor A_h=0$ for
$h\neq j$, their inverses are
$(\pi_{B_i}^R\btensor\id)\,l_{X_p}^{-1}$ and
$(\id\btensor\pi_{B_j}^R)\,r_{X_p}^{-1}$, respectively.

If $p,q$ are not composable, then $X_p\btensor X_q=0$.
When $p=e_i$, the two definitions of $\kappa_{e_i,e_i}$ above agree  by the triangle axiom.
The compatibility of the isomorphisms $\kappa_{p,q}$ with the
associativity constraints is the commutativity of the following
pentagon, in which $\alpha$ denotes the associativity constraint:
\[
\begin{tikzcd}[column sep=large,row sep=large]
(X_p\btensor X_q)\btensor X_r
  \arrow[r,"\kappa_{p,q}\btensor\id_{X_r}"]
  \arrow[d,"\alpha"']
&
X_{pq}\btensor X_r
  \arrow[dr,"\kappa_{pq,r}"]\\
X_p\btensor(X_q\btensor X_r)
  \arrow[r,"\id_{X_p}\btensor\kappa_{q,r}"']
&
X_p\btensor X_{qr}
  \arrow[r,"\kappa_{p,qr}"']
&
X_{pqr}
\end{tikzcd}
\]
for composable paths $p,q,r$, while the compatibility with the unit
constraints is expressed by the identities
$\kappa_{e_i,p}=l_{X_p}(\iota^R_{B_i}\btensor\id_{X_p})$ and
$\kappa_{p,e_j}=r_{X_p}(\id_{X_p}\btensor\iota^R_{B_j})$ used in their
definition. For paths of positive length the pentagon follows from Mac
Lane's coherence theorem. When trivial paths are involved, it follows from
the naturality of the unit constraints and the triangle axiom together
with the definitions of the maps $\kappa$.

Since the noncomposable components vanish, the maps $\kappa_{p,q}$
give a natural isomorphism
$F_{B,X}(M)\btensor F_{B,X}(N)\xrightarrow{\sim}
F_{B,X}(M\btensor N)$, denoted $J^{B,X}_{M,N}$,
while $\bigoplus_iA_i=\mathbf1_R$ gives the unit isomorphism
$J_0^{B,X}$.
The pentagon displayed above, pasted over the decompositions
$M\cong\bigoplus_pM_p\otimes_{\Bbbk}S_p$ and
$N\cong\bigoplus_qN_q\otimes_{\Bbbk}S_q$ and combined with the associativity
constraints of $\mathcal C_Q$ and $\mathcal C_R$, gives the pentagon axiom
for $F_{B,X}$, and the unit identities give the triangle axiom. Thus
$F_{B,X}$ is a strong monoidal functor.

Now let $(F,J,J_0)$ be an arbitrary $\Bbbk$-linear strong monoidal functor.
Let $(B_i)$, $\theta_{e_i}$, and $X_a$ be as in the proof of the preceding
lemma, and set $\theta_a\coloneqq\id_{X_a}$.

Let
$m^Q_{p,q}\colon S_p\btensor S_q\xrightarrow{\sim}S_{pq}$
be the canonical isomorphism, and set
\[
\widehat J^F_{p,q}
\coloneqq
F(m^Q_{p,q})
J_{S_p,S_q}.
\]

For a nontrivial path $p$ and an arrow $a$ with $s(a)=t(p)$, define
recursively
\[
\theta_{pa}
=
\widehat J^F_{p,a}
(\theta_p\btensor\theta_a).
\]
By induction on the length of $q$, using the compatibility of
$\kappa$ with the associativity constraints
together with the associativity and unit coherence axioms for $F$, we obtain, for all composable $p,q$,
\begin{equation}\label{xxeq4-1}
\theta_{pq}\kappa_{p,q}
=
\widehat J^F_{p,q}
(\theta_p\btensor\theta_q).
\end{equation}

Since $F$ is $\Bbbk$-linear, for a finite-dimensional $\Bbbk$-vector space $V$
and $X\in\mathcal C_Q$, the map $v\otimes_{\Bbbk} y\mapsto F(i_v)(y)$, where
$i_v\colon X\to V\otimes_{\Bbbk}X$ sends $x$ to $v\otimes_{\Bbbk} x$, is an isomorphism
$V\otimes_{\Bbbk}F(X)\cong F(V\otimes_{\Bbbk}X)$, natural in $V$ and $X$.
Together with $M\cong\bigoplus_pM_p\otimes_{\Bbbk}S_p$, this gives a canonical
isomorphism
$\bigoplus_pM_p\otimes_{\Bbbk}F(S_p)\cong F(M)$, natural in $M$.

Under this identification, set
$\Theta_M\coloneqq\bigoplus_p(\id_{M_p}\otimes_{\Bbbk}\theta_p)$.
The monoidality of $\Theta$ is the commutativity of the square
\[
\begin{tikzcd}[column sep=large,row sep=large]
F_{B,X}(M)\btensor F_{B,X}(N)
  \arrow[r,"J^{B,X}_{M,N}"]
  \arrow[d,"\Theta_M\btensor\Theta_N"']
&
F_{B,X}(M\btensor N)
  \arrow[d,"\Theta_{M\btensor N}"]\\
F(M)\btensor F(N)
  \arrow[r,"J^F_{M,N}"']
&
F(M\btensor N),
\end{tikzcd}
\]
whose $(p,q)$-component is \eqref{xxeq4-1}, displayed
as
\[
\begin{tikzcd}[column sep=10em,row sep=large]
X_p\btensor X_q
  \arrow[r,"\kappa_{p,q}"]
  \arrow[d,"\theta_p\btensor\theta_q"']
&
X_{pq}
  \arrow[d,"\theta_{pq}"]\\
F(S_p)\btensor F(S_q)
  \arrow[r,"F(m^Q_{p,q})\circ J_{S_p,S_q}"']
&
F(S_{pq}).
\end{tikzcd}
\]
The unit compatibility is the commutativity of
\[
\begin{tikzcd}[column sep=large,row sep=large]
\mathbf1_R
  \arrow[r,"J_0^{B,X}"]
  \arrow[dr,"J_0"']
&
F_{B,X}(\mathbf1_Q)
  \arrow[d,"\Theta_{\mathbf1_Q}"]\\
&
F(\mathbf1_Q),
\end{tikzcd}
\]
which follows from
$J_0\iota_{B_i}^R=F(\iota_i^Q)\theta_{e_i}$ for all $i$, together with
$\bigoplus_iA_i=\mathbf1_R$.
Hence $\Theta\colon F_{B,X}\xRightarrow{\sim}F$ is a monoidal natural
isomorphism.
\end{proof}

\begin{lemma}\label{xxlem4.4}
Let
$F_{B,X}$
and
$F_{B',Y}$
be the functors constructed from the data
$((B_i),(X_a))$
and
$((B_i'),(Y_a))$,
respectively.
If $(B_i)\neq(B_i')$,
then there is no monoidal natural transformation between them.
If $B_i=B_i'$ for every $i$,
then there is a bijection
\[
\operatorname{Hom}_\otimes(F_{B,X},F_{B',Y})
\cong
\prod_{a\in Q_1}
\operatorname{Hom}_{\mathcal C_R}(X_a,Y_a).
\]

\end{lemma}

\begin{proof}
Let $\eta\colon F_{B,X}\Rightarrow F_{B',Y}$
be a monoidal natural transformation.
Write $\kappa^X_{p,q}$ and $\kappa^Y_{p,q}$ for the isomorphisms
$\kappa_{p,q}$ associated with the families $X$ and $Y$, respectively.
Write
$J_0^{B,X}\colon\mathbf1_R\to F_{B,X}(\mathbf1_Q)$ and
$J_0^{B',Y}\colon\mathbf1_R\to F_{B',Y}(\mathbf1_Q)$
for the unit isomorphisms of the two functors.

The unit compatibility condition gives
$\eta_{\mathbf1_Q}\circ J_0^{B,X}=J_0^{B',Y}$.
Applying naturality to the idempotent $E_i^Q$ and using the preceding
equality, we obtain
\[
(J_0^{B',Y})^{-1}
F_{B',Y}(E_i^Q)
J_0^{B',Y}
=
(J_0^{B,X})^{-1}
F_{B,X}(E_i^Q)
J_0^{B,X}.
\]
The two sides are $E_{B_i'}^R$ and $E_{B_i}^R$, respectively, so
$B_i'=B_i$.
Since $B_i=B_i'$ for every $i$, the two functors agree on
$\mathbf1_Q$ and their unit isomorphisms coincide.
The unit compatibility condition therefore gives
$\eta_{\mathbf1_Q}=\id_{F_{B,X}(\mathbf1_Q)}$.
Applying naturality to the inclusion $\iota_i^Q$ gives
$\iota_{B_i}^R\circ\eta_{S_{e_i}}=\iota_{B_i}^R$.
Composing on the left with the direct-sum projection $\pi_{B_i}^R$ gives
$\eta_{S_{e_i}}=\id_{A_i}$.
Set $f_a\coloneqq\eta_{S_a}\colon X_a\to Y_a$ $(a\in Q_1)$.
Using the identifications
$F_{B,X}(S_p)\cong X_p$ and
$F_{B',Y}(S_p)\cong Y_p$,
we regard $\eta_{S_p}$ as a morphism $X_p\to Y_p$.
For a path $p$ and an arrow $a$ with $s(a)=t(p)$, the monoidality of
$\eta$ is expressed by the commutativity of the square
\[
\begin{tikzcd}[column sep=large,row sep=large]
X_p\btensor X_a
  \arrow[r,"\kappa^X_{p,a}"]
  \arrow[d,"\eta_{S_p}\btensor f_a"']
&
X_{pa}
  \arrow[d,"\eta_{S_{pa}}"]\\
Y_p\btensor Y_a
  \arrow[r,"\kappa^Y_{p,a}"']
&
Y_{pa}.
\end{tikzcd}
\]
Since $\kappa^X_{p,a}$ is invertible, $\eta_{S_{pa}}$ is determined by
$\eta_{S_p}$ and $f_a$. Together with $\eta_{S_{e_i}}=\id_{A_i}$,
induction on the length of $p$ shows that every $\eta_{S_p}$ is
determined by the $f_a$.
For $M\in\mathcal C_Q$ and $x\in M_p$, let
$g_x\colon S_p\to M$, $u_p\mapsto x$,
be the corresponding comodule morphism.
Naturality of $\eta$ with respect to $g_x$ shows that,
on the summand $M_p\otimes_{\Bbbk} X_p$, we have
$\eta_M=\id_{M_p}\otimes_{\Bbbk}\eta_{S_p}$.
Hence $\eta$ is uniquely determined by the family
$(f_a)_{a\in Q_1}$.

Conversely, suppose that a family of morphisms
$f_a\colon X_a\to Y_a$ $(a\in Q_1)$ is given.
Set $f_{e_i}=\id_{A_i}$.
For a nontrivial path $p$ and an arrow $a$ with $s(a)=t(p)$, define
recursively
\[
f_{pa}
=
\kappa^Y_{p,a}
(f_p\btensor f_a)
(\kappa^X_{p,a})^{-1}.
\]
Induction on the length of $q$, together with coherence, gives, for all
composable paths $p,q$,
\begin{equation}\label{xxeq4-2}
f_{pq}\kappa^X_{p,q}
=
\kappa^Y_{p,q}
(f_p\btensor f_q).
\end{equation}

For $M\in\mathcal C_Q$, define
$\eta_M\coloneqq\bigoplus_p(\id_{M_p}\otimes_{\Bbbk} f_p)$.
Naturality follows because comodule morphisms preserve the homogeneous
components indexed by paths, while monoidality follows from
\eqref{xxeq4-2}.
The equality
$f_{e_i}=\id_{A_i}$
gives the unit compatibility.
Thus $\eta$ is a monoidal natural transformation.
The two constructions are inverse to each other and are compatible
with composition, giving the stated bijection.
\end{proof}

\begin{proof}[Proof of Theorem~\ref{xxthm4.1}]
The assignment $((B_i),(X_a))\mapsto F_{B,X}$, together with the
assignment on morphisms constructed in
Lemma~\ref{xxlem4.4}, defines a functor
$\Psi$ from the right-hand side of the asserted equivalence to
$\operatorname{Fun}^{\otimes}_{\Bbbk}(\mathcal C_Q,\mathcal C_R)$.
Lemmas~\ref{xxlem4.2}
and~\ref{xxlem4.3}
show that $\Psi$ is essentially surjective, while
Lemma~\ref{xxlem4.4}
shows that $\Psi$ is fully faithful.
Hence $\Psi$ is an equivalence.
\end{proof}

\section{Projection Formulas and Induced Functors on Centers}

Let
$F\colon\mathcal C_Q\to\mathcal C_R$
be a $\Bbbk$-linear strong monoidal functor with right adjoint
$G\colon\mathcal C_R\to\mathcal C_Q$,
and let
$\varepsilon\colon FG\Rightarrow\operatorname{Id}_{\mathcal C_R}$
be the counit.

For
$M\in\mathcal C_Q$
and
$N\in\mathcal C_R$,
the left and right projection formula morphisms are
\begin{align*}
\lproj_{M,N}&\colon
M\btensor G(N)
\longrightarrow
G(F(M)\btensor N),\\
\rproj_{N,M}&\colon
G(N)\btensor M
\longrightarrow
G(N\btensor F(M)).
\end{align*}
They correspond under the adjunction $F\dashv G$, respectively, to
\[
(\id_{F(M)}\btensor\varepsilon_N)
\circ
J_{M,G(N)}^{-1},
\qquad
(\varepsilon_N\btensor\id_{F(M)})
\circ
J_{G(N),M}^{-1}.
\]

We say that the left projection formula holds if every
$\lproj_{M,N}$
is an isomorphism, and that the right projection formula holds if every
$\rproj_{N,M}$
is an isomorphism.

The validity of either projection formula is invariant under monoidal
natural isomorphism of strong monoidal functors. Hence, by
Theorem~\ref{xxthm4.1}, it is enough to consider
the functors $F_{B,X}$ constructed above.

We now specialize to the functor
$F=F_{B,X}\colon\mathcal C_Q\to\mathcal C_R$
of Lemma~\ref{xxlem4.3}, and retain the
decomposition
$R_0=\bigsqcup_{i\in Q_0}B_i$,
the objects $X_p$, and the isomorphisms $\kappa_{p,q}$ from its
construction.

\begin{lemma}\label{xxlem5.1}
The functor $F$ has a right adjoint
$G\colon\mathcal C_R\to\mathcal C_Q$ given by
\begin{equation}\label{xxeq5-1}
G(N)
=
\bigoplus_{p\in\mathcal P(Q)}
\operatorname{Hom}_{\mathcal C_R}(X_p,N)
\otimes_{\Bbbk}S_p.
\end{equation}
In particular, if
$w\in B_i$,
then $G(S_{e_w}^R)\cong S_{e_i}$.

\end{lemma}

\begin{proof}
From $F(M)\cong\bigoplus_{p\in\mathcal P(Q)}M_p\otimes_{\Bbbk}X_p$
and finite dimensionality, there are natural isomorphisms
\[
\operatorname{Hom}_{\mathcal C_R}(F(M),N)
\cong
\bigoplus_{p\in\mathcal P(Q)}
\operatorname{Hom}_{\Bbbk}
\left(
M_p,
\operatorname{Hom}_{\mathcal C_R}(X_p,N)
\right)
\cong
\operatorname{Hom}_{\mathcal C_Q}(M,G(N)).
\]

Thus the functor $G$ defined by
\eqref{xxeq5-1}
is right adjoint to $F$.
The counit $\varepsilon$ is, on each direct summand, the evaluation map
$f\otimes_{\Bbbk} x\longmapsto f(x)$.

If $w\in B_i$, then
$\operatorname{Hom}_{\mathcal C_R}(X_{e_i},S_{e_w}^R)\cong \Bbbk$,
whereas
$\operatorname{Hom}_{\mathcal C_R}(X_{e_h},S_{e_w}^R)=0$
for $h\neq i$.

If $p$ has positive length, then
$s(p)\neq t(p)$, so
$X_p\in\mathcal C_R[B_{s(p)},B_{t(p)}]$
has no simple direct summand indexed by a trivial path.
Indeed, every simple summand $S_r^R$ of $X_p$ satisfies
$s(r)\in B_{s(p)}$ and $t(r)\in B_{t(p)}$. If $r=e_w$ were trivial, then
$w\in B_{s(p)}\cap B_{t(p)}=\varnothing$, a contradiction.
Hence $G(S_{e_w}^R)$ has only one nonzero homogeneous component, namely
the one indexed by $e_i$, and that component is one-dimensional, so
$G(S_{e_w}^R)\cong S_{e_i}$.
\end{proof}

Fix $p\in\mathcal P(Q)$.
If $q=ph$ for some path $h$, then the homogeneous component of
$\lproj_{S_p,N}$ indexed by $q$ is the map
\begin{equation}\label{xxeq5-2}
\operatorname{Hom}_{\mathcal C_R}(X_h,N)
\longrightarrow
\operatorname{Hom}_{\mathcal C_R}(X_{ph},X_p\btensor N),
\qquad
f\longmapsto
(\id_{X_p}\btensor f)
\circ
\kappa_{p,h}^{-1}.
\end{equation}
If $q$ is not of this form, the source is zero.

Similarly, if $q=hp$ for some path $h$, then the homogeneous component of
$\rproj_{N,S_p}$ indexed by $q$ is the map
\begin{equation}\label{xxeq5-3}
\operatorname{Hom}_{\mathcal C_R}(X_h,N)
\longrightarrow
\operatorname{Hom}_{\mathcal C_R}(X_{hp},N\btensor X_p),
\qquad
f\longmapsto
(f\btensor\id_{X_p})
\circ
\kappa_{h,p}^{-1}.
\end{equation}
If $q$ is not of this form, the source is zero.
Both formulas follow directly from the description of the counit as
evaluation.

For each arrow $a\in Q_1$, set
$m_a(r)\coloneqq\dim_{\Bbbk}(X_a)_r$ for $r\in\mathcal P(R)$. Then
\[
X_a
\cong
\bigoplus_{r\in\mathcal P(R)}
\Bbbk^{m_a(r)}
\otimes_{\Bbbk}
S_r^R,
\qquad
\mathcal S_a
\coloneqq
\{
r\in\mathcal P(R)
\mid
m_a(r)>0
\}.
\]

Here and below, for $r\in\mathcal S_a$ we write
$s(r)\in B_{s(a)}$ and $t(r)\in B_{t(a)}$ for the source and target of the
path $r$. Thus $s$ and $t$ also denote the resulting maps from
$\mathcal S_a$ to $B_{s(a)}$ and $B_{t(a)}$.

If $u=rz$ (resp.\ $u=zr$) for some, possibly trivial, path $z$, we call
$r$ a prefix (resp.\ suffix) of $u$. If $z$ is nontrivial, we call $r$ a
proper prefix (resp.\ proper suffix) of $u$.

\begin{theorem}\label{xxthm5.2}
Let $Q,R$ be finite acyclic quivers, let
$F=F_{B,X}\colon\mathcal C_Q\to\mathcal C_R$
be the strong monoidal functor constructed in
Lemma~\ref{xxlem4.3},
and let $G$ be its right adjoint.

\textup{(1)}
The left projection formula holds if and only if, for every arrow
$a\colon i\to j$,
all nonzero multiplicities
$m_a(r)$
are equal to $1$,
the map $t\colon\mathcal S_a\to B_j$ is a bijection, and, for every arrow $b$ of $Q$, every
$r\in\mathcal S_a$, and every
$u\in\mathcal S_b$ with $(a,r)\neq(b,u)$,
the path $r$ is not a prefix of $u$.

\textup{(2)}
The right projection formula holds if and only if, for every arrow
$a\colon i\to j$,
all nonzero multiplicities
$m_a(r)$
are equal to $1$,
the map $s\colon\mathcal S_a\to B_i$ is a bijection, and, for every arrow $b$ of $Q$, every
$r\in\mathcal S_a$, and every
$u\in\mathcal S_b$ with $(a,r)\neq(b,u)$,
the path $r$ is not a suffix of $u$.

\textup{(3)}
Both the left and right projection formulas hold if and only if, for every
arrow $a\colon i\to j$, there exist a bijection
$\sigma_a\colon B_i\xrightarrow{\sim}B_j$ and, for each $v\in B_i$, a path
$\tau_a(v)\colon v\to\sigma_a(v)$ such that
\begin{equation}\label{xxeq5-4}
X_a
\cong
\bigoplus_{v\in B_i}S_{\tau_a(v)}^R,
\end{equation}
and, as $(a,v)$ varies, the paths $\tau_a(v)$ are pairwise distinct, with
none of them a proper prefix or a proper suffix of another.
\end{theorem}

\begin{proof}
For $q\in\mathcal P(Q)$ let
$\operatorname{pr}_q\colon X\to X_q\otimes_{\Bbbk}S_q$ denote the projection
onto the $q$-th homogeneous component. By
Proposition~\ref{xxprop2.1}, every morphism in
$\mathcal C_Q$ preserves homogeneous components. Hence, for all
$M\in\mathcal C_Q$ and $N\in\mathcal C_R$, the square
\[
\begin{tikzcd}[column sep=6em,row sep=large]
M\btensor G(N)
  \arrow[r,"\lproj_{M,N}"]
  \arrow[d,"\operatorname{pr}_q"']
&
G(F(M)\btensor N)
  \arrow[d,"\operatorname{pr}_q"]\\
(M\btensor G(N))_q\otimes_{\Bbbk}S_q
  \arrow[r,"(\lproj_{M,N})_q\otimes_{\Bbbk}\id_{S_q}"']
&
\bigl(G(F(M)\btensor N)\bigr)_q\otimes_{\Bbbk}S_q
\end{tikzcd}
\]
commutes, where $(\lproj_{M,N})_q$ denotes the induced map on
$q$-components. Consequently, if $\lproj_{M,N}$ is an isomorphism, then so
is $(\lproj_{M,N})_q$ for every $q$.

We first prove the necessity of the multiplicity and target
conditions in \textup{(1)}.
Fix
$a\colon i\to j$
and
$w\in B_j$.
By
Lemma~\ref{xxlem5.1},
\[
G(S_{e_w}^R)\cong S_{e_j},
\qquad
S_a\btensor G(S_{e_w}^R)\cong S_a.
\]

On the other hand, $F(S_a)=X_a$, and the homogeneous component indexed by
$a$ in $G(X_a\btensor S_{e_w}^R)$ is
${\operatorname{Hom}_{\mathcal C_R}(X_a,X_a\btensor S_{e_w}^R)}$.

By the formula for tensor products of simple objects,
\[
\dim_{\Bbbk}
\operatorname{Hom}_{\mathcal C_R}(X_a,X_a\btensor S_{e_w}^R)
=
\sum_{\substack{r\in\mathcal S_a\\t(r)=w}}
m_a(r)^2.
\]
If the left projection formula holds, then by the component square above
(with $q=a$) the map is an isomorphism between two one-dimensional spaces,
and the right-hand side must be equal to $1$.
Since it is a sum of squares of nonnegative integers, there is exactly one
$r\in\mathcal S_a$
such that $t(r)=w$, and this $r$ satisfies $m_a(r)=1$.
As $w\in B_j$ was arbitrary, all nonzero multiplicities are $1$, and
$t\colon\mathcal S_a\to B_j$ is a bijection.

We next prove the necessity of the prefix condition.
Suppose that there exist an arrow $b$ of $Q$ and paths $r,u,z$ such that
\[
r\in\mathcal S_a,
\qquad
u\in\mathcal S_b,
\qquad
u=rz,
\qquad
(a,r)\neq(b,u).
\]
Then $S_u^R$ occurs as a direct summand of both
$X_b$ and $X_a\btensor S_z^R$. Choose a split inclusion
$\iota_b\colon S_u^R\to X_b$ and a split projection
$p_b\colon X_b\to S_u^R$ of the first summand, and let
$\iota_a\colon S_u^R\to X_a\btensor S_z^R$ be the inclusion of the second.
Then $g\coloneqq\iota_a\circ p_b\neq0$:
\[
\begin{tikzcd}[column sep=large,row sep=large]
&
S_u^R
  \arrow[dl,hook,bend left=20,"\iota_b"']
  \arrow[dr,hook,"\iota_a"]
&\\
X_b
  \arrow[ur,bend left=20,"p_b"]
  \arrow[rr,dashed,"g=\iota_a p_b\neq0"']
& &
X_a\btensor S_z^R.
\end{tikzcd}
\]
Consequently,
\[
\bigl(G(X_a\btensor S_z^R)\bigr)_b\neq0,
\qquad
\bigl(S_a\btensor G(S_z^R)\bigr)_b=0.
\]
Indeed, $b$ is an arrow, so $\ell(b)=1$. If $b=ah$, then $\ell(h)=0$,
hence $h=e_{t(a)}$ and $a=b$.
Thus, if $b\neq a$, then $b$ cannot be written as $ah$ for any path $h$.
If $b=a$ and $b=ah$, then necessarily $h=e_j$.
Since $r\neq u$ in this case, $z$ is nontrivial, and therefore
$\operatorname{Hom}_{\mathcal C_R}(X_{e_j},S_z^R)=0$.
This contradicts the assumption that
$\lproj_{S_a,S_z^R}$
is an isomorphism.
Hence $r$ is not a prefix of $u$.

We now prove sufficiency.
Assume the conditions in \textup{(1)}.
Let $p=a_1\cdots a_m\colon i\to j$ be a nontrivial path, and let $w\in B_j$.
By the bijectivity of the target maps, there are unique paths
$r_\ell\in\mathcal S_{a_\ell}$ $(1\leq\ell\leq m)$ such that
\[
t(r_m)=w,
\qquad
t(r_\ell)=s(r_{\ell+1})
\quad (1\leq\ell<m).
\]
Define $\tau_p^{-}(w)\coloneqq r_1\cdots r_m$.
For a trivial path $e_i$ and $w\in B_i$, set
$\tau_{e_i}^{-}(w)\coloneqq e_w$.

Using the formula for tensor products of simple objects, induction on
$\ell(p)$ then gives
\begin{equation}\label{xxeq5-5}
X_p
\cong
\bigoplus_{w\in B_{t(p)}}
S_{\tau_p^{-}(w)}^R.
\end{equation}

By construction, for composable paths $q_1,q_2$ and
$w\in B_{t(q_2)}$,
\begin{equation}\label{xxeq5-6}
\tau_{q_1q_2}^-(w)
=
\tau_{q_1}^-\bigl(s(\tau_{q_2}^-(w))\bigr)\,
\tau_{q_2}^-(w).
\end{equation}

Suppose $\tau_p^-(w)$ is a prefix of $\tau_q^-(w')$.
We show that $p$ is a prefix of $q$.
Since an arrow $a\colon i\to j$ has $i\neq j$, the sets $B_i$ and $B_j$
are disjoint. As $X_a\in\mathcal C_R[B_i,B_j]$, every path in
$\mathcal S_a$ is nontrivial.
Hence the cases in which $p$ or $q$ is trivial follow immediately from
$\tau_{e_i}^-(w)=e_w$ and the disjointness of the sets $B_i$.
We may therefore assume that $p$ and $q$ are nontrivial, and write
$p=a_1\cdots a_m$ and $q=b_1\cdots b_n$.
By construction, we may write
$\tau_p^-(w)=r_1\cdots r_m$ and
$\tau_q^-(w')=u_1\cdots u_n$, where
$r_\alpha\in\mathcal S_{a_\alpha}$ $(1\leq\alpha\leq m)$ and
$u_\beta\in\mathcal S_{b_\beta}$ $(1\leq\beta\leq n)$.

Since $r_1$ is a prefix of $\tau_p^-(w)$, which is itself a prefix of
$\tau_q^-(w')$, and $u_1$ is also a prefix of $\tau_q^-(w')$,
one of $r_1$ and $u_1$ is a prefix of the other.
If $r_1=u_1$ but $a_1\neq b_1$, then $(a_1,r_1)\neq(b_1,u_1)$ although
$r_1$ is a prefix of $u_1$, contradicting the prefix condition.
Hence $r_1=u_1$ forces $a_1=b_1$. Combined with the case $r_1\neq u_1$,
this gives $a_1=b_1$ and $r_1=u_1$.
Cancelling this common prefix and repeating the argument, we obtain
$m\leq n$ and $a_\alpha=b_\alpha$ for $1\leq\alpha\leq m$. Hence $p$ is a prefix of $q$.
If $q=ph$, then after cancelling the common prefix corresponding to
$p$, the remaining path is $\tau_h^-(w')$.

In particular, if $\tau_p^-(w)=\tau_q^-(w')$, then $p=q$ and
$w=w'$. Thus, for every $r\in\mathcal P(R)$, there is at most one
$q\in\mathcal P(Q)$ such that
$r=\tau_q^-(t(r))$.
Hence
\begin{equation}\label{xxeq5-7}
G(S_r^R)
\cong
\begin{cases}
S_q,
&
r=\tau_q^{-}(t(r))
\text{ for some }q\in\mathcal P(Q),
\\
0,
&
\text{otherwise}.
\end{cases}
\end{equation}

Fix $p\in\mathcal P(Q)$ and $r\in\mathcal P(R)$.
If $s(r)\notin B_{t(p)}$, then $X_p\btensor S_r^R=0$.
In this case, by \eqref{xxeq5-7}, if
$G(S_r^R)\cong S_h\neq0$, then
$s(h)\neq t(p)$, and hence $S_p\btensor S_h=0$.
Thus the source and target of
$\lproj_{S_p,S_r^R}$
are both zero.

Suppose now that $s(r)\in B_{t(p)}$.
Let $u=\tau_p^{-}(s(r))$.
By \eqref{xxeq5-5}, $X_p\btensor S_r^R\cong S_{ur}^R$.
By \eqref{xxeq5-7}, the preceding prefix argument, and
\eqref{xxeq5-6}, $ur=\tau_q^{-}(t(r))$
for some $q\in\mathcal P(Q)$ if and only if
$G(S_r^R)\cong S_h$ for some $h$ with $s(h)=t(p)$,
in which case $q=ph$.

Thus
\[
G(X_p\btensor S_r^R)
\cong
\begin{cases}
S_{ph},
&
G(S_r^R)\cong S_h,
\quad
s(h)=t(p),
\\
0,
&
\text{otherwise},
\end{cases}
\]
which is isomorphic to $S_p\btensor G(S_r^R)$.

In the nonzero case, the source and target of
\eqref{xxeq5-2}
are both one-dimensional.
Choose a direct-sum decomposition $X_h=Y\oplus Y'$ with
$Y\cong S_r^R$, and let $f\colon X_h\to S_r^R$ be the projection
onto $Y$ followed by a fixed isomorphism
$Y\xrightarrow{\sim}S_r^R$. Let $\iota_Y\colon Y\to X_h$ be the
inclusion. Set $Z\coloneqq\kappa_{p,h}(X_p\btensor Y)$, a direct summand of
$X_{ph}$ that $\kappa_{p,h}^{-1}$ carries isomorphically onto
$X_p\btensor Y$.
Since $s(r)\in B_{t(p)}$, we have
$X_p\btensor Y\cong X_p\btensor S_r^R\cong S_{ur}^R\neq0$, and the
composite
\[
\begin{tikzcd}[column sep=large,row sep=large]
X_{ph}
  \arrow[r,"\kappa_{p,h}^{-1}"]
&
X_p\btensor X_h
  \arrow[r,"\id_{X_p}\btensor f"]
&
X_p\btensor S_r^R
  \arrow[r,"\sim"]
&
S_{ur}^R\neq0\\
Z
  \arrow[u,hook]
&
X_p\btensor Y
  \arrow[l,"\sim"']
  \arrow[u,hook,"\id_{X_p}\btensor\iota_Y"']
  \arrow[ur,"\sim"']
&
&
\end{tikzcd}
\]
is nonzero. Hence $(\id_{X_p}\btensor f)\circ\kappa_{p,h}^{-1}\neq0$, so the map in
\eqref{xxeq5-2} is nonzero and therefore, between
one-dimensional spaces, an isomorphism.
Thus $\lproj_{S_p,S_r^R}$ is an isomorphism.
Passing to finite direct sums and using $\Bbbk$-linearity, we conclude that every
$\lproj_{M,N}$
is an isomorphism.

The proof of \textup{(2)} is analogous to that of \textup{(1)}:
using \eqref{xxeq5-3} in place of
\eqref{xxeq5-2}, the same argument applies with the
source condition in place of the target condition and the suffix condition
in place of the prefix condition.

Finally, suppose that both the left and right projection formulas hold.
Then, for each arrow $a\colon i\to j$, the source and target maps are bijections.
Hence one can write uniquely
\[
\mathcal S_a
=\{\tau_a(v)\mid v\in B_i\},
\qquad
s(\tau_a(v))=v,
\qquad
t(\tau_a(v))=\sigma_a(v),
\]
where $\sigma_a\colon B_i\xrightarrow{\sim}B_j$.

This gives
\eqref{xxeq5-4}.
To see that the paths $\tau_a(v)$ are pairwise distinct, suppose that
$\tau_a(v)=\tau_b(u)$ for some $(a,v)\neq(b,u)$. Taking
$r=u=\tau_a(v)$ in condition \textup{(1)} contradicts the prefix
condition, since $r$ is then a prefix of $u$. The prefix and suffix
conditions likewise exclude proper prefix and proper suffix relations
among these paths. Thus the paths $\tau_a(v)$ satisfy the condition in
\textup{(3)}, and the conditions in \textup{(3)} are necessary.

Conversely, suppose that the conditions in \textup{(3)} hold.
Then all nonzero multiplicities are equal to $1$, and the maps
$s\colon\mathcal S_a\to B_i$ and
$t\colon\mathcal S_a\to B_j$
are bijections for every arrow $a\colon i\to j$.
The condition on the paths $\tau_a(v)$ in \textup{(3)} gives the prefix and
suffix conditions in \textup{(1)} and \textup{(2)}.
Hence both projection formulas hold.
\end{proof}

\begin{examples}\label{xxexample5.3}
\textup{(1)} Let $Q$ consist of one arrow $a\colon i\to j$, and let $R$ consist
of one arrow $c\colon x\to y$. Take $B_i=\{x\}$, $B_j=\{y\}$, and
$X_a=S_c^R\oplus S_c^R$.
The source and target conditions hold, and the prefix and suffix conditions
are vacuous, but $m_a(c)=2$. Hence neither projection formula holds.

\textup{(2)} Let $Q$ have arrows
$a\colon i\to j$ and $b\colon i\to\ell$, and let
$R$ be the quiver
\[
x\xrightarrow{\,c\,}y\xrightarrow{\,d\,}z.
\]
Set $B_i=\{x\}$, $B_j=\{y\}$, $B_\ell=\{z\}$,
$X_a=S_c^R$, and $X_b=S_{cd}^R$.
All nonzero multiplicities are $1$, and the source and target maps are
bijective.
Since $c$ is a proper prefix of $cd$, the left projection formula does not
hold, whereas the right projection formula holds because neither $c$ nor
$cd$ is a proper suffix of the other.

\textup{(3)} Let $Q$ have two parallel arrows $a,b\colon i\to j$, and let
$R$ consist of one arrow $c\colon x\to y$. Take $B_i=\{x\}$,
$B_j=\{y\}$, and $X_a=X_b=S_c^R$.
All nonzero multiplicities are $1$, and the source and target maps are
bijective. Moreover, the condition in
Theorem~\ref{xxthm5.2}\textup{(1)} would
remain satisfied
if ``prefix'' were replaced by ``proper prefix'', since $c$ is not a
proper prefix of itself. Nevertheless neither projection formula holds:
for the pair $(a,c)\neq(b,c)$ the path $r=u=c$ lies in
$\mathcal S_a\cap\mathcal S_b$, so $r$ is a non-proper prefix (and also a
non-proper suffix) of $u$, and the conditions in
Theorem~\ref{xxthm5.2}\textup{(1)} and
\textup{(2)} fail.
\end{examples}

\begin{corollary}\label{xxcor5.4}
Suppose that both the left and right projection formulas hold, and let
$\sigma_a$ and $\tau_a$ be as in
Theorem~\ref{xxthm5.2}\textup{(3)}.
Then $G$ induces a functor
$\widetilde G\colon Z_{\mathrm{lax}}(\mathcal C_R)\to
Z_{\mathrm{lax}}(\mathcal C_Q)$.
Let $\Phi_Q$ and $\Phi_R$ be the equivalences of
Theorem~\ref{xxthm3.2} for $Q$ and $R$,
respectively.
Then there is a natural isomorphism
$\Phi_Q\circ\widetilde G\cong T\circ\Phi_R$, where the functor
$T\colon\operatorname{Rep}^{\mathrm{fd}}_{\Bbbk}(R)\to
\operatorname{Rep}^{\mathrm{fd}}_{\Bbbk}(Q)$
is given on objects by
\begin{align}
T(V)_i
&=
\bigoplus_{v\in B_i}V_v,
\label{xxeq5-8}\\
\varphi_a^{T(V)}\big|_{V_v}
&=
\varphi^V_{\tau_a(v)}\colon
V_v
\longrightarrow
V_{\sigma_a(v)}
\qquad
(a\colon i\to j),
\label{xxeq5-9}
\end{align}
and on morphisms by $T(f)_i=\bigoplus_{v\in B_i}f_v$.

Moreover, the restriction of $\widetilde G$ to the Drinfeld centers is a braided lax
monoidal functor $Z(\mathcal C_R)\to Z(\mathcal C_Q)$.
Under the equivalences above, this braided lax monoidal functor corresponds
to the restriction of $T$ to
$\operatorname{Rep}^{\mathrm{fd}}_{\Bbbk,\mathrm{iso}}(R)$, whose lax monoidal
structure at a vertex $i$ is given by the projection
$T(V)_i\otimes_{\Bbbk}T(W)_i\to T(V\otimes W)_i$ onto the diagonal direct
summands, with unit morphism $\lambda\mapsto(\lambda)_{v\in B_i}$.
The functor $T$ also admits a separable Frobenius monoidal structure whose
lax monoidal part is given by the same formulas.
\end{corollary}

\begin{proof}
Let $(V,\gamma)\in Z_{\mathrm{lax}}(\mathcal C_R)$.
For $M\in\mathcal C_Q$, the map $\widetilde\gamma_M$ is the
right vertical arrow in the following commutative diagram:
\begin{equation}\label{xxeq5-10}
\begin{tikzcd}[column sep=large,row sep=large]
G(V\btensor F(M))
  \arrow[d,"G(\gamma_{F(M)})"]
&
G(V)\btensor M
  \arrow[l,"\rproj_{V,M}"]
  \arrow[d,"\widetilde\gamma_M"']\\
G(F(M)\btensor V)
  \arrow[r,"(\lproj_{M,V})^{-1}"']
&
M\btensor G(V)
\end{tikzcd}
\end{equation}
By \cite[Proposition~3.22]{flake2026projection}, the
inverses of the projection formula morphisms make $G$ a
$\mathcal C_Q$-bimodule functor
$\mathcal C_R\to\mathcal C_Q$, where
$\mathcal C_R$ is regarded as a $\mathcal C_Q$-bimodule category via $F$.
Naturality of $\widetilde\gamma$ follows from the naturality of
$\lproj$, $\rproj$, and $\gamma$.
For the unit condition, the unit coherence of the bimodule-functor
structure on $G$, together with the unit condition for $\gamma$, gives
the commutativity of
\[
\begin{tikzcd}[column sep=large,row sep=large]
G(V)\btensor\mathbf1_Q
  \arrow[r,"\widetilde\gamma_{\mathbf1_Q}"]
  \arrow[dr,"r_{G(V)}"']
&
\mathbf1_Q\btensor G(V)
  \arrow[d,"l_{G(V)}"]\\
&
G(V),
\end{tikzcd}
\]
that is,
$\widetilde\gamma_{\mathbf1_Q}=l_{G(V)}^{-1}\circ r_{G(V)}$.

For the compatibility condition, naturality of $\gamma$ with respect to
$J_{M,N}$ identifies the occurrence of $\gamma_{F(M\btensor N)}$ in the
definition of $\widetilde\gamma_{M\btensor N}$ with the corresponding map for
$F(M)\btensor F(N)$.
The compatibility condition for $\gamma$ expresses the latter as the
composite of the two successive half-braidings.
Applying the left and right module-functor coherence to the projection
formula morphisms, followed by the bimodule-functor compatibility, gives
\[
\begin{tikzcd}[column sep=large,row sep=large]
&
G(V)\btensor M\btensor N
  \arrow[dl,"\widetilde\gamma_M\btensor\id_N"']
  \arrow[dr,"\widetilde\gamma_{M\btensor N}"]
&\\
M\btensor G(V)\btensor N
  \arrow[rr,"\id_M\btensor\widetilde\gamma_N"']
& &
M\btensor N\btensor G(V),
\end{tikzcd}
\]
that is,
$\widetilde\gamma_{M\btensor N}
=(\id_M\btensor\widetilde\gamma_N)
\circ(\widetilde\gamma_M\btensor\id_N)$.
Hence $(G(V),\widetilde\gamma)$ is an object of
$Z_{\mathrm{lax}}(\mathcal C_Q)$.

If $f\colon(V,\gamma)\to(W,\delta)$ is a morphism in the lax center,
naturality of the projection formula morphisms and the defining condition for $f$
show that $G(f)$ intertwines $\widetilde\gamma$ and
$\widetilde\delta$.
Thus the construction defines the functor $\widetilde G$.
As shown in the first part of the proof of
Theorem~\ref{xxthm3.2},
$V$ is supported only on the trivial paths of $R$.
Hence, by
\eqref{xxeq5-1},
$G(V)_{e_i}=\bigoplus_{v\in B_i}V_{e_v}$,
and all homogeneous components indexed by paths of positive length vanish.
This gives
\eqref{xxeq5-8}.

Take $M=S_a$, where $a\colon i\to j$.
The homogeneous component of $G(V\btensor X_a)$ indexed by $a$ is
$\operatorname{Hom}_{\mathcal C_R}(X_a,V\btensor X_a)$.
The isomorphism
$\rproj_{V,S_a}$
in \eqref{xxeq5-10} identifies the summand $V_{e_v}$
with the summand corresponding to $\tau_a(v)$.
On this summand,
$G(\gamma_{X_a})$
is given by the linear map
$\varphi^V_{\tau_a(v)}$.
The isomorphism
$(\lproj_{S_a,V})^{-1}$
in \eqref{xxeq5-10} then identifies the
result with $V_{\sigma_a(v)}$.
This gives
\eqref{xxeq5-9}.
The formula on morphisms follows similarly by taking direct sums over the vertices.

Under the identifications
$G(V)_{e_i}=\bigoplus_{v\in B_i}V_{e_v}$, the canonical lax monoidal map
$G(V)\btensor G(W)\to G(V\btensor W)$ is the projection onto the diagonal
direct summands, since
$S_{e_v}^R\btensor S_{e_w}^R=0$ whenever $v\neq w$, and its unit morphism
is the diagonal map $\lambda\mapsto(\lambda)_{v\in B_i}$.
By \cite[Theorem~4.10]{flake2026projection}, the lax
monoidal structure on the induced functor between Drinfeld centers is
the one inherited from $G$.
Thus, under $\Phi_Q$ and $\Phi_R$, it is the lax monoidal structure on the
restriction of $T$ described in the statement.
We extend this structure to a separable Frobenius monoidal structure on
$T$. Denote its component at a vertex $i$ by
\[
\mu^T_{V,W,i}\colon
\bigoplus_{v,w\in B_i}
V_v\otimes_{\Bbbk}W_w
\longrightarrow
\bigoplus_{v\in B_i}
V_v\otimes_{\Bbbk}W_v.
\]
Let $\delta^T_{V,W,i}$ be the inclusion in the opposite direction, and let
the oplax unit morphism be given at each vertex $i$ by
$\bigoplus_{v\in B_i}\Bbbk\to \Bbbk$,
$(\lambda_v)_{v\in B_i}\mapsto\sum_{v\in B_i}\lambda_v$.

The maps $\mu^T_{V,W}$ and $\delta^T_{V,W}$ are natural in $V$ and
$W$.
Since each $\sigma_a$ is a bijection, the map associated with an arrow
$a\colon i\to j$ sends the summand indexed by $(v,w)$ to the summand indexed by
$(\sigma_a(v),\sigma_a(w))$, and
$v=w$ if and only if $\sigma_a(v)=\sigma_a(w)$.
Thus the diagonal and off-diagonal summands are preserved, so
$\mu^T_{V,W}$ and $\delta^T_{V,W}$ are morphisms of
$Q$-representations.
For the same reason, the lax and oplax unit morphisms
are morphisms of $Q$-representations.
The unitality and counitality axioms follow
directly from the definitions of these two unit morphisms.
Associativity, coassociativity, and the two Frobenius compatibility
relations can be verified on the relevant direct-sum components: both sides
are the identity when the corresponding indices coincide and are zero
otherwise.
Hence $T$ is a Frobenius monoidal functor.
Moreover,
$\mu^T_{V,W}\circ\delta^T_{V,W}=\id_{T(V\otimes W)}$,
which at each vertex $i$ is the commutativity of
\[
\begin{tikzcd}[column sep=large,row sep=large]
\bigoplus_{v\in B_i}V_v\otimes_{\Bbbk}W_v
  \arrow[r,"\delta^T_{V,W,i}"]
  \arrow[dr,"\id"']
&
\bigoplus_{v,w\in B_i}V_v\otimes_{\Bbbk}W_w
  \arrow[d,"\mu^T_{V,W,i}"]\\
&
\bigoplus_{v\in B_i}V_v\otimes_{\Bbbk}W_v.
\end{tikzcd}
\]
so $T$ is a separable Frobenius monoidal functor
\cite[Section~2]{MS2010}.

Finally, if $V\in\operatorname{Rep}^{\mathrm{fd}}_{\Bbbk,\mathrm{iso}}(R)$,
then every arrow in $R$, and hence every path
$\tau_a(v)$,
acts invertibly.
The map in
\eqref{xxeq5-9}
is a direct sum of these invertible maps arranged according to the bijection
$\sigma_a$,
and is therefore invertible.
Thus $T$ restricts as claimed.
\end{proof}

\begin{example}\label{xxexample5.5}
The following example shows that both projection formulas may hold although
$F$ is not an equivalence.
Let $Q$ consist of one arrow $a\colon i\to j$, and let $R$
consist of two arrows $c\colon x\to y$ and $d\colon u\to v$. Take
$B_i=\{x,u\}$, $B_j=\{y,v\}$, and $X_a=S_c^R\oplus S_d^R$.
The maps $s\colon\mathcal S_a\to B_i$ and $t\colon\mathcal S_a\to B_j$
are bijections, and the conditions in
Theorem~\ref{xxthm5.2}\textup{(3)} hold with
$\sigma_a(x)=y$, $\sigma_a(u)=v$, $\tau_a(x)=c$, and $\tau_a(u)=d$. Hence
both projection formulas hold. However,
$F(S_{e_i})\cong S_{e_x}^R\oplus S_{e_u}^R$ is not simple, so $F$ is not
an equivalence. Under the equivalence of
Corollary~\ref{xxcor5.4}, the induced functor $T$
satisfies $T(V)_i=V_x\oplus V_u$, $T(V)_j=V_y\oplus V_v$, and
$\varphi_a^{T(V)}=\varphi_c^V\oplus\varphi_d^V$.
\end{example}

\subsection*{Funding} 
Jingheng Zhou was partially supported by NSFC (No. 12301053).

\subsection*{General AI Use}
The authors have used AI to improve the readability and language of the manuscript. The authors confirm that the content of the manuscript was not generated by AI and that they take full responsibility for the content of the publication.

\providecommand{\bysame}{\leavevmode\hbox to3em{\hrulefill}\thinspace}
\providecommand{\MR}{\relax\ifhmode\unskip\space\fi MR }
% \MRhref is called by the amsart/book/proc definition of \MR.
\providecommand{\MRhref}[2]{%
  \href{http://www.ams.org/mathscinet-getitem?mr=#1}{#2}
}
\providecommand{\href}[2]{#2}

\end{document}